\documentclass [hidelinks,12pt,openbib]{article}

\usepackage{cancel}
\usepackage{xcolor}
\usepackage[hidelinks]{hyperref}
\usepackage[normalem]{ulem}
 
\usepackage{mathtools}
 \usepackage{tikz}
 \usepackage{cancel}
\usepackage{fixmath}
\usepackage{dsfont}
\usepackage{caption}
\usepackage{pslatex}
\usepackage{bm}
\usepackage{graphicx}
\usepackage{amsmath,amssymb,amsfonts,amsthm}
\usepackage{xcolor}
\usepackage{mathtools}
\usepackage{dsfont}
\usepackage{cite}
\usepackage{ulem}
\usepackage{fancybox,fancyvrb}
\usepackage{framed} 

\usepackage{MnSymbol,wasysym}
\usepackage[export]{adjustbox}

\usepackage{enumitem}
\usepackage{color}
\usepackage{fancyhdr}
\usepackage{esint}
\usepackage{mathtools}
\usepackage[margin=1in]{geometry}
\usepackage{mdframed}

\def\XXint#1#2#3{{\setbox0=\hbox{$#1{#2#3}{\int}$}
     \vcenter{\hbox{$#2#3$}}\kern-.5\wd0}}
 \def\begfrwhite{\begin{mdframed}[
	backgroundcolor=white!15!white,
	linecolor=black,
	linewidth=0.5pt,
	align=center,
	userdefinedwidth=4.99in]}

\def\E1{\text{E}_1}

\def\L1loc{{L^1_{\rm loc}}}

\def\bc{\begin{center}}
\def\ec{\end{center}}

\def\begfr{\begin{mdframed}[
	backgroundcolor=gray!15!white,
	linecolor=blue,
	linewidth=1pt,
	align=center,
	userdefinedwidth=4.99in]}
\def\begfrwhite{\begin{mdframed}[
	backgroundcolor=white!15!white,
	linecolor=black,
	linewidth=0.5pt,
	align=center,
	userdefinedwidth=4.99in]}

\def\begfrblue{\begin{mdframed}[
	backgroundcolor=blue!15!white,
	linecolor=blue,
	linewidth=1pt,
	align=center,
	userdefinedwidth=4.99in]}

\def\endfr{\end{mdframed}}
\def\endfrred{\end{mdframed}}
\def\begfrred{\begin{mdframed}[
	backgroundcolor=red!15!white,
	linecolor=blue,
	linewidth=1pt,
	align=center,
	userdefinedwidth=4.99in]}
\def\endfrblue{\end{mdframed}}

\def\begeq{\begin{equation}}
\def\endeq{\end{equation}}

\def\bbr{\begin{bmatrix*}[r]}
\def\ebr{\end{bmatrix*}}
\def\bb{\begin{bmatrix*}[r]}
\def\eb{\end{bmatrix*}}
\def\bbc{\begin{bmatrix}}

\def\ebc{\end{bmatrix}}

\definecolor{darkmagenta}{rgb}{0.55, 0.0, 0.55}

\definecolor{darkgreen}{rgb}{0.0,0.6,0.0}

\def\0u{\underline{0}}

\def\begeq{\begin{equation}} 
\def\endeq{\end{equation}}

\def\bcenter{\begin{center}}
\def\ecenter{\end{center}}
\def\beq{\begin{equation}}
\def\eeq{\end{equation}}
\def\bmatr{\begin{bmatrix*}[r]}
\def\ematr{\end{bmatrix*}}
\def\bmatc{\begin{bmatrix}}
\def\ematc{\end{bmatrix}}

\usepackage{amsmath}               
  {
      \theoremstyle{plain}
      
  }
  {
      \theoremstyle{plain}
      
  }
   {
      \theoremstyle{plain}
      
  }
   {
      \theoremstyle{plain}
      
  }
  {
      \theoremstyle{plain}
      
  }
    {
      \theoremstyle{plain}
      
  }
\newcommand{\vertiii}[1]{{\left\vert\kern-0.25ex\left\vert\kern-0.25ex\left\vert #1 
    \right\vert\kern-0.25ex\right\vert\kern-0.25ex\right\vert}}
 \usepackage{amsthm}
\newtheorem*{proposition*}{Proposition}
\newtheorem*{corollary*}{Corollary}
\newtheorem*{definition*}{Definition}
\newtheorem*{example*}{Example}
\newtheorem*{lemma*}{Lemma}
\newtheorem*{theorem*}{Theorem}
\newtheorem*{observation*}{Observation}
\newtheorem*{remark*}{Remark}
\usepackage{sectsty}
\sectionfont{\fontsize{15}{10}\selectfont}

\newtheorem*{exmp*}{Example}
\usepackage{stmaryrd}
\usepackage{wrapfig}

\usepackage{xurl} 

\begin{document}

\begin{center}
{\bf \Large Variations on the two-child problem } 
\vskip 15pt
Christoph B\"orgers$^1$ and Samer Nour Eddine$^{2}$
\end{center}

\noindent
Departments of $^1$Mathematics and $^2$Psychology, Tufts University, Medford, MA

\vskip 20pt

\begin{quote}
{\small {\bf Abstract.} } Mr.\ Smith has two children. Given that at least one of them is a boy, how likely is it that Mr.\ Smith has two boys? It's a very standard puzzle
in elementary books on probability theory. Whoever asks
you this question hopes that you will answer ``$\frac{1}{2}$", in which case they can say triumphantly ``Oh no, the answer is $\frac{1}{3}$". 
This is called the {\em two-child puzzle}. 
Some authors have discussed a striking variation, which we'll call the {\em Adam puzzle}. Again, Mr. Smith has two children.  Given that one
of them is a boy named Adam, how likely is it that Mr.\ Smith has two boys? Astonishingly, now the answer is $\frac{1}{2}$, at least approximately.
(The exact answer depends a bit on precise assumptions.) 
We give pictorial explanations of both puzzles. 
We  then point out that the answers usually given rely on a tacit assumption
about how the information that one of Mr.\ Smith's two children is a boy, or one of them is a boy named Adam, is obtained.  We give 
examples showing 
 that the answers may be different with different assumptions. 
We conclude with a discussion of why the 
Adam puzzle is so confusing to most people. 
\end{quote}

\section{Introduction} 

This article is about puzzles concerning conditional probabilities. They have been written about many times before, even in the 
{\em Guardian} newspaper \cite{guardian}. 
The first and best-known is the {\em two-child puzzle}, 
which appeared in Martin Gardner's October 1959  ``Mathematical games column" in {\em Scientific American} \cite{Gardner}. It 
is now found in most introductory textbooks on probability (for instance \cite[Example 4.25]{Grinstead_Snell}).

\begin{quote}
{\em 
{\bf Two-child puzzle:}
Mr.\ Smith has two children. At least one of them is a boy. What is the probability that Mr.\ Smith has two boys?
}
\end{quote}
\noindent
The answer  is $\frac{1}{3}$, not $\frac{1}{2}$ as you might have thought. 

The following variation was discussed for instance in \cite{DAgostini} and  \cite{Drunkards_Walk}. We will call it the {\em Adam puzzle} in this article. 

\begin{quote}
{\em 
{\bf  Adam puzzle:} 
Mr.\ Smith has two children. One of them is a boy named Adam. What is the probability that  Mr.\ Smith has two boys?
}
\end{quote}

\noindent
To everybody's utter surprise, now the answer is $\frac{1}{2}$, approximately or exactly depending on precise assumptions about how parents choose names. 
Naming the boy made a difference, even though we knew in advance that the boy has {\em some} name! 

There is no trick here.
One-third of two-child
families with a son have two boys, and approximately one-half of two-child families with an Adam have two boys. You can verify 
this empirically, or --- more conveniently --- by computer simulation. We'll give you all the required code, five lines in total,  later in this article.
We will explain these answers, and draw pictures that may help you grasp why they are the answers. 

Then we will make two additional points. First, in deriving these answers, a tacit
assumption is made about how you obtained your information about Mr.\ Smith's children. The assumption is that
you got your information from an all-knowing oracle, or that it's just given in the problem.  More real-life-compatible assumptions
sometimes give answers that are different, though not necessarily less surprising.

Second, we discuss why the Adam puzzle is so confusing to many people. We believe that the root of the confusion may lie in 
reading the sentence ``Given that Mr.\ Smith has a boy named Adam, 
the probability that he has two boys is $\frac{1}{2}$" as  saying ``If Mr.\ Smith has a boy named Adam, the probability that he has two boys
is $\frac{1}{2}$."

Silly though these puzzles are, they   are worth thinking about. Conditional probabilities are important. Misconceptions about them 
have resulted in bad medical decisions
 \cite{Gigerenzer_et_al}, in politicians unknowingly or knowingly spreading misinformation (see the story about Rudy Giuliani in \cite{Gigerenzer_et_al}), 
  in innocent people being sent to jail 
 \cite{Prosecutors_Fallacy}, perhaps also in murderers going free (see the story about O.\ J.  Simpson in  \cite{OJ}). Among all the things an undergraduate
 mathematics student should learn, ``Be careful with conditional probabilities" may be the most important and consequential in the real world. 
 
Throughout this article, we make two false assumptions:  (1) A child is as likely a boy as a girl. (2) All children are either boys or girls. Unfortunately we have been unable to find equally memorable ways of phrasing these puzzles avoiding those assumptions.

\section{The two-child puzzle}

The puzzle was stated in the introduction. 
We will give several explanations of the textbook answer $\frac{1}{3}$. 

\subsection{Verbal explanation.}
Of the four {\em a priori} equally likely possibilities boy-boy, boy-girl, girl-boy, and girl-girl (this indicates the sex of each child, in birth order), one was eliminated when  we told you that Mr.\ Smith has a boy: girl-girl is no longer possible. This leaves three equally likely possibilities: boy-boy, boy-girl, and girl-boy. In only one of these three cases does Mr.\ Smith
have a second boy.

\subsection{Pictorial explanation.} 

 Picture the four kinds of families, boy-boy, boy-girl, girl-boy, and girl-girl, as boxes of equal size. One box, the girl-girl box, is ruled
out when you learn that at least one of Mr.\ Smith's two children is a boy. Among the three boxes that remain, only one correponds
to boy-boy families.

\begin{center}
\includegraphics[scale=0.75]{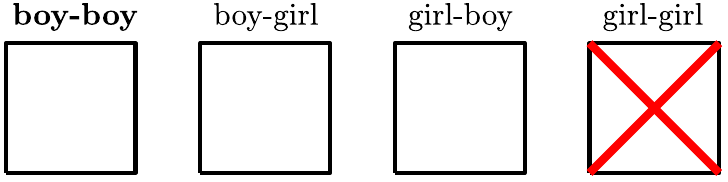}

Figure 1. The boxes indicate different ``types" of families.
\end{center}

\subsection{Two lines of code.} 

If you are like us and nothing convinces you quite as much as a computer simulation, then run these two lines of {\scshape Matlab}  code:
\begin{verbatim}
     N=10^7; X=randi([0,1],2,N); S=sum(X); 
     a=length(find(S==2)); b=length(find(S>=1)); a/b
\end{verbatim}

\noindent
We simulate $N=10^7$ families with two children.  $X(k,i)=1$ if in the $i$-th family, the $k$-th child is a boy, and 
$X(k,i)=0$ otherwise, $k \in \{1,2\}$, $1 \leq i \leq N$.
The number of boys in the $i$-th family is $S(i)$.
The number of families with two boys is $a$, the number of families with at least one boy is $b$. When we
ran the code we found that a fraction of 0.33326... of families with at least one boy had two boys.

\section{The Adam puzzle} 

The puzzle was stated in the introduction. 
We will give several explanations of the answer $\frac{1}{2}$. Under the precise assumptions that we make here, which we'll state shortly, it is actually 
slightly smaller than $\frac{1}{2}$. 

\subsection{Verbal explanation.} You are asking ``What is the sex of the child not named Adam in the Smith family?" The chances
should be even!  

\subsection{Algebraic explanation.} Unless otherwise noted, we always 
assume that a newborn boy will be called ``Adam" with probability $p>0$, unless he already has an older brother by that name. 
Notice that therefore in a two-boy family, the first boy's name is Adam with probability $p$, and the second boy's name is Adam with probability 
$(1-p)p$, 

Consequently the fraction of two-boy families that have an Adam is, altogether, 
\begin{equation}
\label{eq:2pmp2}
p + (1-p) p = 2p-p^2.
\end{equation}
Among families with an older boy and a younger girl, of course a fraction $p$ have an Adam. The same holds for families with an older
girl and a younger boy. Therefore altogether, the fraction of families with an Adam which have two boys is 
\begin{equation}
\label{eq:fraction_of_boys}
\frac{2p-p^2}{2p-p^2+p+p} = \frac{1-\frac{p}{2}}{2 - \frac{p}{2}}
\end{equation}
This expression is slightly smaller than $\frac{1}{2}$. The easy way to see that is to write 
$$
\frac{1-\frac{p}{2}}{2 - \frac{p}{2}} = \frac{2-\frac{p}{2}}{2 - \frac{p}{2}} - \frac{1}{2 - \frac{p}{2}} = 1  - \frac{2}{4 - p}.
$$
As $p$ rises from $0$ to $1$, this expression decreases from $\frac{1}{2}$ to $\frac{1}{3}$. It is not surprising that we should 
get $\frac{1}{3}$ for $p=1$: If any boy is called Adam unless he has an older brother by that name already, then having a boy
and having and Adam are the same thing, so the Adam puzzle becomes the two-child puzzle. 

\subsection{Does it matter that nobody gives the same name to two of their children?} 
 Not at all. To see this, assume  that even if a family already has a boy named Adam,
 if they have another boy, the probability that that boy, too, will be named Adam is still $p$. Now suppose that {\em after the fact}, we required 
 all the families with two Adams to re-name one of their sons, to enforce the ``no two children with the same name in one family" rule. This would not change the set of families who have an Adam 
 at all, and therefore would not change the answer to our puzzle. 
 
 \subsection{Different assumptions on families' child-naming habits may result in slightly different answers.} 
 Suppose that there are $n$ possible names for boys, and for a family's first boy, each has probability $p=\frac{1}{n}$ of being chosen, but
 then for the second boy (if there is one), one of the remaining $n-1$ names is chosen, each having probability $\frac{1}{n-1}$. The probability that 
 a two-boy family has an Adam is now 
 $$
 \frac{1}{n} + \left( 1 - \frac{1}{n} \right) \frac{1}{n-1} = p + (1-p) \frac{p}{1-p} = 2p,
 $$
 and consequently the probability that a family with an Adam has two boys is now 
 $$
 \frac{2p}{2p+p+p} = \frac{1}{2}
 $$
 exactly. 
 
 We can do a similar calculation assuming, more generally and more realistically, that the $n$ available names are chosen with different probabilities 
 $p_1,\ldots,p_n$, and that in a two-boy family who chose the $k$-th name for their first boy, the probabilities for the remaining names will 
 be boosted by  the factor $\frac{1}{1-p_k}$ for a second boy. The answer to the Adam puzzle now slightly depends on $p_1,\ldots, p_n$, but is close to $\frac{1}{2}$ if 
 all $p_k$ are small. 
 
For the remainder of this article, we return to our simple initial assumption: Any newborn boy is called Adam with probability $p$, unless he already has an older brother
 by that name. 
 
\subsection{Pictorial explanation.} 

When you are told that Mr.\ Smith has a boy named Adam, you are told that his family
is not just in the left three boxes, but specifically in a small area within those boxes, indicated in blue here: 
\begin{center}
\includegraphics[scale=0.75]{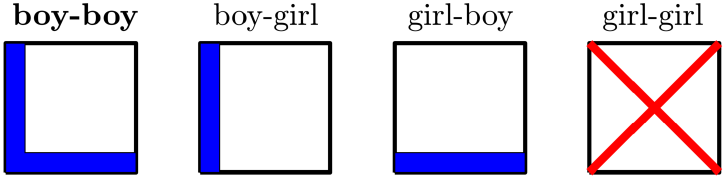}

Figure 2. Families with boys named Adam are indicated in blue.
\end{center}

\noindent
Think of the boxes as being of size $1 \times 1$, and the blue stripes as being of width $p$. 
Then the expression in  (\ref{eq:2pmp2}) represents the size of the 
blue area in the boy-boy box. It is almost as large as the blue in the other two boxes combined. This is why the chance of having 
two boys is  nearly $\frac{1}{2}$.

\subsection{Three lines of code.}
This time it takes three lines of {\scshape Matlab} code:

\begin{verbatim}
     N=10^9; X=randi([0,1],2,N); S=sum(X);
     U=max(0,randi(50,2,N)-49); V=U.*X; A=sum(V); 
     a=length(find(A>0 & S>1)); b=length(find(A>0)); a/b
\end{verbatim}

The variables $N$, $X$, and $S$ have the same meaning as before, with $N$ indicating the number of two-child families considered, 
the binary matrix $X$ indicating which child in which family is a boy, and the vector $S$ indicating how many boys each family has. 
Notice that we are now running a billion samples. 
That took a little while on our laptop. You can try $N=10^7$ and already get very similar results. 

We assume here that families are allowed to name two boys Adam. As we have explained, that 
assumption does not change anything, as long as we assume that the probability that Adam is chosen 
as the name of a boy is always $p$.
Entries in $U$ are either $1$ or $0$, and they are $1$ with probability 
$ p = \frac{1}{50}$. If $U(k,i)=1$, we assume that the $k$-th child in the $i$-th family will be
named Adam {\em if} that child is a boy. Then 
$V$ is obtained from $U$ by setting entries corresponding to girls to zero, so $V(k,i)=1$ if
the $k$-th child in the $i$-th family is a boy named Adam. 
The number of Adams in the $i$-th family is $A(i)$. It can be 0, 1, or 2. Then 
$a$ is the number of families who have an Adam and another boy, and $b$ is the number of families who have 
an Adam. When we ran this code we found that a fraction 0.49747... of 
two-child families with an Adam had two boys.  Formula (\ref{eq:fraction_of_boys}) with $p=\frac{1}{50}$ yields 0.4975.

\subsection{But really, why does the reasoning change when you mention the boy's name?} 

In the original two-child puzzle, the reasoning is based on the fact that 
the possibilities that were not ruled out (boy-boy, boy-girl, and girl-boy) are still
equally likely. This is no longer the case once you mention the boy's name, as shown by Figure 2.
Families with two boys are more likely to have an Adam than families with just one boy, just as having two raffle tickets
makes you more likely to win than having only one. Therefore (strictly speaking by 
Bayes' formula), 
when you are told that Mr.\ Smith has an Adam, that hints at two boys.

 \section{It matters how you learned about Mr.\ Smith's children} 
 
 In both puzzles, we know that Mr.\ Smith has two children, one of them a boy, or one of them a boy named Adam, because the puzzle said so. 
 Let's think about more real-life ways of obtaining that information. For simplicity, we {\em are} going to assume that we already know Mr.\ Smith has two
 children. Our question is how we learn about the sex, and possibly the name, of one of them.

 \subsection{Parent meeting at school.}  Consider the following variation.
 
 \begin{quote}
{\em 
Mr. Smith  has two children. 
You are at the parent meeting at school. The principal says ``Would those families who have boys please go to room 333 to hear about the 
athletics programs for boys?" You see Mr.\ Smith get up and head towards room 333. He has a boy! 
What is the probability 
that Mr. Smith’s other child is also a boy? 
}
\end{quote}

The answer is clearly $\frac{1}{3}$. (Well, that's clear because you read our earlier discussion, otherwise most of us would have said $\frac{1}{2}$.) This is the two-child puzzle. In effect, the principal elicited the answer to the question ``Do you have a boy?" By heading to 
room 333, Mr.\ Smith revealed that for him, the answer is yes.  (We are ignoring the  possibility that Mr.\ Smith might have a 2-year-old boy who isn't yet going to school, or a 27-year-old boy who is no longer going to school.) 

\subsection{Accidental encounter with Mr.\ Smith and his son.} 
The following variation was discussed in \cite{bar_hillel_and_falk}: 

\begin{quote}
{\em 
Mr.\ Smith has two children. You meet him walking along the street with
a young boy whom he proudly introduces as his son. What is the probability 
that Mr. Smith’s other child is also a boy? 
}
\end{quote}

Most people find the answer to this question obvious. It is exactly $\frac{1}{2}$, as explained pictorially in Figure 3.\footnote{Tacitly, we assume here that Mr.\ Smith 
selected a child for his walk at random, without any bias. If for instance Mr.\ Smith would never take a girl on a walk, then 
seeing him with a boy merely tells us that he has at least one boy; the probability that he has a second boy is still $\frac{1}{3}$ in that case.}

\begin{center}
\includegraphics[scale=0.75]{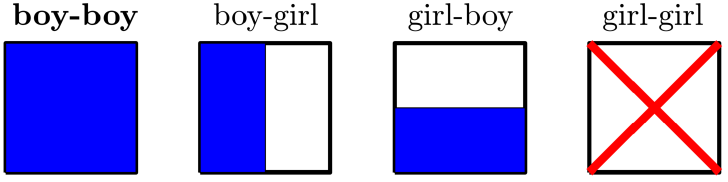}

Figure 3. Instances when Mr.\ Smith's randomly selected child is a boy  in blue.
\end{center}

Intuitively, why does the chance encounter of Mr.\ Smith with his son make a second boy more likely? The answer is that 
encountering Mr.\ Smith with a son becomes more likely if he has two boys than if he only has one. Therefore encountering him with a son hints at two boys (strictly 
speaking by Bayes' formula). 

\subsection{Accidental encounter with Mr.\ Smith and his son Adam.} In the previous puzzle, let's change one word: Mr.\ Smith introduces his son by name. 

\begin{quote}
{\em 
Mr.\ Smith has two children. You meet him walking along the street with
a young boy whom he proudly introduces as his son Adam. What is the probability 
that Mr. Smith’s other child is also a boy? 
}
\end{quote}

Again assume that 
Mr.\ Smith selected a random child to accompany him. If he has two boys, the probability of selecting an Adam for the walk was 
$$
\frac{1}{2} p + \frac{1}{2} p (1-p).
$$
The factors of $\frac{1}{2}$ represent the likelihood of choosing the first boy or the second boy. The first boy is named Adam with probability $p$. 
The second boy is named Adam with probability $p(1-p)$. If Mr.\ Smith has an older boy and a younger girl, the probability of selecting an Adam for the walk is
$$
\frac{1}{2} p. 
$$
The same is true if he has an older girl and a younger boy. Therefore if you see him walk with a boy named Adam, the probability that his is a two-boy family is
$$
\frac{\frac{1}{2} p + \frac{1}{2} p (1-p)}{\frac{1}{2} p + \frac{1}{2} p (1-p) + \frac{1}{2} p + \frac{1}{2} p } = \frac{1- \frac{p}{2}}{2 - \frac{p}{2}}.
$$
It's precisely the answer derived before for the Adam puzzle. The probability, which was {\em exactly} $\frac{1}{2}$ before, is now very slightly smaller than 
$\frac{1}{2}$. We should probably pay little attention to that change, since it depends on the exact assumptions we make on how parents select names
for their children, as explained earlier.

\subsection{Military draft.} Finally let's think about the following puzzle. 

\begin{quote}
{\em Mr.\ Smith has two children. One day the government announces that all families who have at least one boy must send a boy, of their choice if they have more than one, to register for the draft. You see  Mr.\ Smith arrive at the registration office in company of a boy, who 
introduces himself as ``Adam Smith". What is the probability that Mr.\ Smith has two boys?} 
\end{quote} 

It appears at first sight that you just learned that Mr.\ Smith has a boy named Adam, no less and no more, and that you are asking ``What is the sex of the non-Adam
at home?" So the answer should be $\frac{1}{2}$, shouldn't it? 

Let's assume that the families with two boys choose the boy who gets sent in for draft registration by flipping a fair coin, paying no attention to their names.
 The probability that a two-boy family will send in an Adam is therefore 
$$
\frac{1}{2} p + \frac{1}{2} p (1-p).
$$
The probability that a family with an older boy and a younger girl will send in an Adam is 
$p$,
and the same holds for a family with a family with an older girl and a younger boy. 
Therefore the probability that a family sending in an Adam has two boys is 
$$
\frac{\frac{1}{2} p + \frac{1}{2} p (1-p)}{\frac{1}{2} p + \frac{1}{2} p (1-p)+p + p}  =  \frac{1-\frac{p}{2}}{3 - \frac{p}{2}}
$$
As $p$ rises from $0$ to $1$, this decreases from $\frac{1}{3}$ to $\frac{1}{5}$ (!) --- rather different from the $\frac{1}{2}$ that we might have expected. 

\section{Why is the Adam puzzle so counter-intuitive?}

\subsection{A counter-argument.} Perhaps the following counter-argument 
captures why many people are confused by the Adam puzzle. It certainly captures why 
{\em we} were first confused by it.

\begin{quote}
{\em 
If Mr.\ Smith has a boy named Adam,  the probability that he has two boys is about $\frac{1}{2}$, The same must
then be true, of course, 
if he has a boy named Sam, or a boy named Victor, or a boy with any other first name. 

\vskip 5pt
\noindent
If Mr.\ Smith has at least one
boy, he has a boy with {\em some} first name, and therefore the probability that he has two boys ought to be 
about $\frac{1}{2}$. 

\vskip 5pt
\noindent
How, then, can the solutions to the Adam puzzle and the original two-child puzzle be different?
}
\end{quote}

This argument is based on the following misunderstanding.

\subsection{Reading conditional probabilities as implications.}  
The confusion likely originates from the way we talk about conditional probabilities. We might say something like this:
\begin{quote}
{\em If Mr.\ Smith has a boy named Adam, the probability that he has two boys is  $\frac{1}{2}$.}
\end{quote}
Let generally $E$ and $F$ be events, and $P(E|F)$ the conditional probability that $E$ occurs, given that $F$ occurs. 
One might easily think that 
\begin{equation}
\label{eq:cp} 
P(E|F) = q
\end{equation}
(where $q$ is some number between $0$ and $1$) could be re-stated as 
\begin{equation}
\label{eq:ip}
F ~~\Rightarrow~~ P(E)=q. 
\end{equation}
However, while (\ref{eq:cp}) makes sense, (\ref{eq:ip}) does not. 

To explain, we remind 
the reader that events are {\em sets} in probability theory. We always think of an underlying random experiment. The set
of all possible outcomes of this experiment is commonly denoted by $\Omega$, and is called the sample space.  In the Adam puzzle, 
the random experiment is ``choose a random two-child family". The sample space is the set of all two-child families. 
An event $E$ is a subset of $\Omega$. You can also think of it as the event that the
outcome of the random experiment lies in $E$. For instance, if $E$ is the event that  the randomly chosen two-child
family has two boys, then strictly speaking we should think of $E$ as the {\em set} of all two-child families with two boys.

By definition, $P(E|F)$ is the size of $E \cap F$, divided by the size of $F$. In rigorous probability theory,  we scale so that 
$\Omega$ has size 1, and we say ``measure" instead of
``size".  Here we can simply take ``size" to mean the number of elements. So $P(E|F)=q$ means that the size of the set $E \cap F$ is $q$ times the size of $F$. 

On the other hand, what should (\ref{eq:ip}) mean?  The arrow ``$\Rightarrow$" ought to stand between statements, not sets. Remembering that $F$ is a set, we should replace the ``$F$"
in (\ref{eq:ip}) by a statement, perhaps like this: Denoting by $\omega \in \Omega$ the outcome of our random experiment, 
\begin{equation}
\label{eq:ip_2}
\omega \in F ~~ \Rightarrow ~~ P(E) = q.
\end{equation}
What's strange about (\ref{eq:ip_2})  is that $P(E)=q$ has nothing to do with any one particular $\omega$. The meaning
of $P(E)=q$ is that a fraction $q$ of all possible $\omega \in \Omega$ also belong to $E$. We don't mean that. What we really mean is that a fraction $q$ of the
times when $F$ occurs, $E$ also occurs. In other
words, we mean  (\ref{eq:cp}). 

\subsection{A misconception following from reading conditional expectations as implications.} 

Suppose $F_1,F_2,\ldots,F_n$ and $E$ are events, and at least one of the $F_k$ must occur. Suppose that for all 
$k$, 
$$
P(E|F_k) = q
$$
for the same number $q$ between 0 and 1. Can we conclude that $P(E) = q$? The counter-argument that we presented earlier
is based on the assumption that we can. 

If we (mis)read $P(E|F_k)=q$ as $F_k \Rightarrow P(E) = q$, then of course that would be correct reasoning --- at least 
one of the $F_k$ most ``hold" (thinking of $F_k$ as a statement rather than an event), and therefore $P(E)$ must be $q$.
However, reading $P(E|F_k)=q$ as $F_k \Rightarrow P(E) = q$ is wrong, and in fact the conclusion is wrong. 

As an example, suppose you throw a dart at a square, and assume that the dart is equally likely to land anywhere 
in the square,  so its landing location is uniformly distributed in the square. 

\begin{center}
\includegraphics[scale=.8]{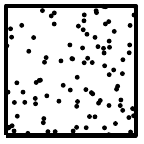}

Figure 4. The square at which we throw darts, with 100 examples of points where the dart  might land.

\end{center}

\noindent
We let $E$ be the event that the dart lands in the right upper quadrant of the square. Clearly $P(E)=\frac{1}{4}$. 

\begin{center}
\includegraphics[scale=.8]{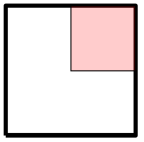}

Figure 5. The dart will land in the red area with probability $\frac{1}{4}$. 

\end{center}

\noindent
We now let $F_1$ and $F_2$ be the events indicated in blue and grey in the following figure. 
Clearly at least one of $F_1$ and $F_2$ occurs, and 
$$
P(E|F_1) = P(E|F_2) = \frac{1}{3}, ~~~\mbox{but} ~~~ P(E) = \frac{1}{4}.
$$

\begin{center}
\includegraphics[scale=.8]{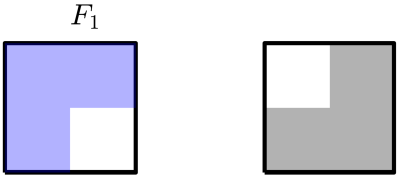}

Figure 6. The ``event $F_1$" is the event that the dart lands in the blue, and the \\ ``event $F_2$" is the event that it lands in the 
grey.

\end{center}

\subsection{The missing condition is disjointness.} 
If  $F_1,\ldots,F_n$ are {\em disjoint} events (whenever one occurs, none of the others occur), and one of them must always occur, 
and $P(E|F_k)=q$ for all $k$, then indeed $P(E) = q$. This is a fact in elementary probability theory. We won't spell out its proof here, 
but Figure 7 shows an example.

\begin{center}
\includegraphics[scale=.8]{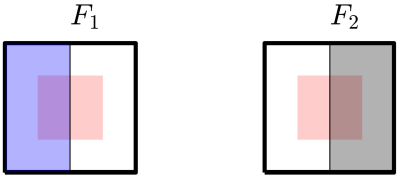}

Figure 7. Here $F_1$ and $F_2$ are disjoint, one of them must occur, and the target has been centered. Now $P(E|F_1) = P(E|F_2) = \frac{1}{4}$, and indeed $P(E) = \frac{1}{4}$.
\end{center}

\subsection{Disjointness and the Adam puzzle.} 

For the counter-argument that we gave earlier to be correct, ``He has a boy named Adam", 
``He has a boy named Sam", and ``He has a boy named Victor" would have to be disjoint events. They are not, of course. 
He can have
a boy named Adam, and another named Sam.
If you assumed that they were disjoint events, that would bring you to the following (wrong) mental picture:

\begin{center}
\includegraphics[scale=0.75]{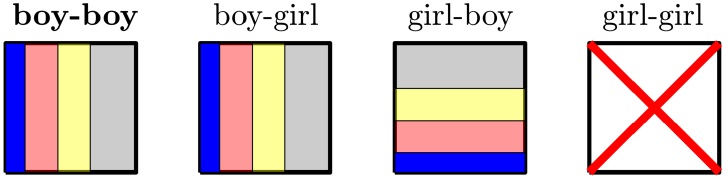}

Figure 8. Different colors correspond to different names. Blue means ``Adam". 

It's the wrong picture for the Adam puzzle! 
\end{center}
\noindent
It would be the right picture if the puzzle were: 

\begin{quote}
{\em 
Mr.\ Smith has two children. At least one of them is a boy. His first-born boy [not his first-born {\em child}, but his first-born {\em boy}] is named Adam. 
What is the probability that  Mr.\ Smith has  two boys?
}
\end{quote}

\noindent
``His first-born boy is named Adam" and ``His first-born boy is named Sam" are in fact disjoint (we disregard
double names, obviously), and Figure 8 is the correct picture for this much less interesting version of the puzzle.
The answer is now  $\frac{1}{3}$, as for the original two-child puzzle.

For the Adam puzzle, where all you are told is that
Mr.\ Smith has two children, and he has a boy named Adam, the correct picture is drawn in Figure 2, and again for several names (four
in total to keep the figure uncluttered) in Figure 9.

\begin{center}
\includegraphics[scale=0.75]{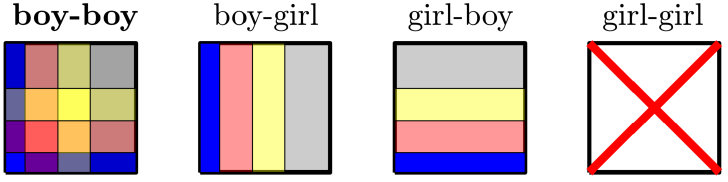}

Figure 9. Different colors correspond to different names. Blue means ``Adam".

This is the right picture for the Adam puzzle.
\end{center}

\vskip 10pt
\noindent
{\bf Acknowledgments.} 
We thank   Ed Larsen, Tilman B\"orgers, Robert Krasny  (all three at the University of Michigan), and Benjamin Borgers for helpful, clarifying, and 
entertaining
conversations about these puzzles.

\end{document}